\documentclass[11pt]{amsart}
\usepackage{amsmath,amsthm,amscd,amssymb, color}
\usepackage{latexsym}
\usepackage[colorlinks, citecolor = blue,backref]{hyperref}
\usepackage[alphabetic]{amsrefs}
\usepackage[capitalize, nameinlink]{cleveref}
\usepackage{enumerate}
\usepackage[margin=1.333in]{geometry}
\usepackage{pgfplots}
\usepackage{booktabs}
\usepackage[T1]{fontenc}

\newcommand{\Q}{\mathbb Q}
\newcommand{\C}{\mathbb C}
\newcommand{\Z}{\mathbb Z}

\renewcommand{\phi}{\varphi}
\newcommand{\eps}{\varepsilon}

\newcommand{\calO}{\mathcal O}

\newcommand{\bmx}{\left( \begin{matrix}}
\newcommand{\emx}{\end{matrix} \right)}

\newcommand{\new}{\mathrm{new}}

\newcommand{\Sq}{\mathrm{Sq}}
\newcommand{\Orb}{\mathrm{Orb}}
\newcommand{\Gal}{\mathrm{Gal}}

\usepackage{bbm}

\newtheorem{lem}{Lemma}
\numberwithin{lem}{section}

\newtheorem{thm}[lem]{Theorem}

\newtheorem{conj}[lem]{Conjecture}
\newtheorem{question}[lem]{Question}
\crefname{question}{Question}{Questions}

\theoremstyle{remark}

\theoremstyle{definition}

\numberwithin{equation}{section}

\begin{document}

\title{An on-average Maeda-type conjecture in the level aspect}
\author{Kimball Martin}
\email{kimball.martin@ou.edu}
\address{Department of Mathematics, University of Oklahoma, Norman, OK 73019 USA}

\date{\today}

\maketitle

\begin{abstract}
We present a conjecture on the average number of Galois orbits of newforms
when fixing the weight and varying the level.
This conjecture implies, for instance, that the central $L$-values
(resp.\ $L$-derivatives) are nonzero for $100\%$ of even weight 
prime level newforms with root number $+1$ (resp.\ $-1$).
\end{abstract}


\section{Introduction}


Let $S_k(N)$ (resp.\ $S^\new_k(N)$) be the space of weight $k$
elliptic cusp (resp.\ new) forms of level $\Gamma_0(N)$.
For $S$ a Hecke-stable subspace of $S_k(N)$ such that the set of
normalized eigenforms is closed under the action of $\Gal(\bar \Q/\Q)$, 
denote by $\Orb(S)$ the set of Galois orbits of normalized eigenforms in $S$.

Maeda's conjecture states that for any even $k$ such that $S_k(1) \ne 0$,
$\# \Orb(S_k(1)) = 1$, and in fact that $T_2$ (or any $T_p$) 
acts irreducibly (over $\Q$) on $S_k(1)$ with Galois group of type $S_n$.
More generally,  let $\Sq_r$ denote the set of squarefree positive integers
with exactly $r$ prime factors.  Then for $N \in \Sq_r$, 
$S_k(N)$ has $2^r$ Atkin--Lehner eigenspaces.  
A generalization of Maeda's conjecture \cite{tsaknias} states that, for fixed
$N \in \Sq_r$, one has
$\# \Orb(S_k^\new(N)) = 2^r$ for all $k$ sufficiently large.  
It follows from trace formula methods that $\# \Orb(S^\new_k(N)) \ge 2^r$ 
for all sufficiently large $k$, and this can be made effective \cite{me:dim}.

On the other hand, if we fix $k=2$ and vary $N \in \Sq_r$, 
one expects a strict inequality $\# \Orb(S_2(N)) > 2^r$ infinitely often
due to the existence of sufficiently many elliptic curves
(or abelian surfaces of GL(2) type) of squarefree level.
However, we predict the following on-average analogue of Maeda's
conjecture in the level aspect.

\begin{conj} \label{gal-conj}   Let $k \ge 2$ be even.
Then the average number of Galois orbits of $S_k^\new(N)$ over all $N \in \Sq_r$
is $2^r$, i.e., 
\begin{equation} \label{eq:conj}
\lim_{X \to \infty} \frac {\sum_{N \in \Sq_r(X)} \# \Orb(S_k^\new(N))}
{\# \Sq_r(X) } = 2^r,
\end{equation}
where $\Sq_r(X) = \{ N \in \Sq_r : N \le X \}$.
In fact, for a fixed prime $p$, $T_p$ acts irreducibly on
each of the $2^r$ Atkin--Lehner eigenspaces of $S_k^\new(N)$ for
$100\%$ of $N \in \Sq_r$ coprime to $p$. 
\end{conj}

In particular, this asserts that each Atkin--Lehner eigenspace is spanned
by a single Galois orbit 100\% of the time, analogous to the usual Maeda conjecture.

One application of Maeda's conjecture is to non-vanishing
of central $L$-values for newforms of full level 
\cite{conrey-farmer} (see also \cite[Corollary 2]{kohnen-zagier}).
The above conjecture similarly has applications to non-vanishing
$L$-values and derivatives.  We recall that if a newform $f \in S_k(N)$
has root number $-1$, then the central $L$-value $L(k,f)$ is forced to
vanish by the functional equation.

\begin{thm} \label{thm1}
Let $k \ge 2$ be even, and assume \cref{gal-conj}.  Let $\mathcal F_k'$
denote the collection of newforms in $\bigcup_{N \text{ prime}} S_k(N)$, 
partially ordered by level. Then

\begin{enumerate}
\item for $100\%$ of $f \in \mathcal F_k'$ with root number $+1$,
we have $L(k, f) \ne 0$; and

\item for $100\%$ of $f \in \mathcal F_k'$  with root number $-1$,
we have $L'(k, f) \ne 0$.
\end{enumerate}
\end{thm}

Let us expand a bit on the context of these statements.

First, Maeda's conjecture for full level was formulated on the basis of computational
data and recent work verifies Maeda's conjecture for $k \le 14000$
(\cite{ghitza-mcandrew}, \cite{bengoechea}).
See \cite{tsaknias} and \cite{DPT} for evidence for the analogue for squarefree level
as well as a possible generalization to non-squarefree level.
Our conjecture is based on both heuristics in \cref{sec:heuristics} and data in
\cref{sec:data}.

One could also attempt to generalize \cref{gal-conj} to non-squarefree levels.
In the case of non-squarefree level, the naive expectation for the number of 
non-CM Galois orbits should
be given by the number of possible local representation types at ramified places,
which varies with both the primes $p | N$ as well as $v_p(N)$.  
We do not consider this here, in part because it would be difficult to generate a 
convincing amount of data for non-squarefree level and partly because even
the correct generalization of Maeda's conjecture to non-squarefree level is not clear
(see \cite{DPT}).  Similarly, we do not consider nontrivial nebentypus or odd
weights, but it is reasonable to expect an analogue of \cref{gal-conj} in these settings
as well.

We also discuss the possibility of stronger statements
than \cref{gal-conj}.  For instance, assuming a conjecture of Roberts \cite{roberts} on the
finitude of rational newforms in large weights, we suggest that for large weights
one may have an exact equality $\# \Orb(S_k^\new(N)) = 2^r$ for all but finitely many
$N \in \Sq_r$ (\cref{conj:finite}), which would be an exact analogue of Tsaknias'
generalized Maeda conjecture in the level aspect.  We also briefly discuss how often 
Hecke polynomials are irreducible or Galois groups are of type 
$S_n$ in \cref{sec:irred,sec:galois}.  In addition, we raise some questions
related to \cite{roberts}, \cite{murty} and \cite{KSW} based on our data
(\cref{q:fin,q:ratl}).

Lipnowski and Schaeffer \cite{LS} formulated a conjecture in
a similar vein as \cref{gal-conj}, that
the rational Hecke modules of each Atkin--Lehner eigenspace are
asymptotically simple.  Restricting to $N$ prime, this means that
for $\eps = \pm 1$ the maximal size of an orbit 
in the Atkin--Lehner $\eps$-eigenspace $S_k^{\new,\eps}(N)$ should be 
asymptotic to the dimension of $S_k^{\new,\eps}(N)$, 
i.e.,
\begin{equation} \label{eq:LS}
\lim_{N \to \infty} \frac {\max \{ \# \calO: \calO \in \Orb(S_k^{\new,\eps}(N)) \} }{\dim S_k^{\new,\eps}(N)} = 1.
\end{equation}
While neither \eqref{eq:conj} nor \eqref{eq:LS} imply the other,
\eqref{eq:conj} is morally stronger in that it is suggestive of 
\eqref{eq:LS} but not conversely.  Namely, \eqref{eq:conj}
implies \eqref{eq:LS} holds for a density 1 subsequence of primes, but
\eqref{eq:LS} does not imply any bounds on the number of 
Galois orbits, even restricting to some density 1 subset of primes.

We also remark that \cref{gal-conj} and \eqref{eq:LS} assert very strong
statements about the growth of degrees of rationality fields of newforms.  
For instance it is known that, as the level grows, the proportion of newforms 
with rationality fields of bounded degree tends to 0 \cite{binder}.  Moreover, there exist 
sequences of newforms with rationality fields of degrees growing at least
logarithmically in the level (e.g., \cite{BM}, \cite{BPR}).  \cref{gal-conj} and
\eqref{eq:LS} would imply that the growth of rationality fields is generically linear
in the level.

Regarding \cref{thm1}, the conclusion is expected from minimalist type conjectures
(see \cite{brumer} for $k=2$) and the Katz--Sarnak philosophy.  The
non-vanishing of central $L$-values $L(k, f)$ and $L$-derivatives $L'(k,f)$
is already known for a positive proportion of
such $f$ by \cite{vanderkam}, \cite{kowalski-michel} and
\cite{iwaniec-sarnak}.  
In \cite{iwaniec-sarnak}, a connection is made between large proportions of 
non-vanishing of central values (via a more refined density conjecture)
and Landau--Siegel zeroes.  This suggests
a connection between the average number of Galois orbits and 
Landau--Siegel zeroes.

The proof of \cref{thm1} is an immediate application of the above non-vanishing
results.  The obstruction to extending this to squarefree level is that
one needs the existence of non-vanishing $L$-values or $L$-derivatives in 100\% 
of the Atkin--Lehner eigenspaces.  While this may very well be accessible by current 
analytic methods, to our knowledge it has not been considered.
One similarly has applications of \cref{gal-conj} to questions such as non-vanishing
of twisted $L$-values.

In \cite{me:zeroes}, we give another application of \cref{gal-conj} to zeroes of
automorphic forms on definite quaternion algebras.

\subsection*{Acknowledgements}
We thank Ariel Pacetti, 
 David Roberts, Gabor Wiese, and especially Bartosz Naskr\k{e}cki  for helpful discussions and feedback.
Some of the computing for this project was performed at the OU Supercomputing Center for Education \& Research (OSCER) at the University of Oklahoma (OU). 
The author was supported by a grant from the 
Simons Foundation/SFARI (512927, KM).


\section{Heuristics}
\label{sec:heuristics}

\subsection{Random polynomials} \label{sec:randpoly}
Let $k \ge 2$ be even, $N \in \Sq_r$ and 
$n = \dim S_k^\new(N) \sim \frac{(k-1)\phi(N)}{12}$.   
Let $\eps$ denote a sign pattern for $N$, i.e., a collection
$(\eps_p)_{p | N}$ such that $\eps_p = \pm 1$ for each $p | N$.
For $p | N$, let $W_p$ be the Atkin--Lehner involution at $p$ on $S_k^\new(N)$.
Denote by $S_k^{\new,\eps}(N) = \{ f \in S_k^\new(N) : W_p f = \eps_p f \text{ for } p | N \}$ the Atkin--Lehner eigenspace of newforms associated to $\eps$.
By \cite{me:dim}, we know that $\dim S_k^{\new, \eps}(N) \sim \frac n{2^r}$
(asymptotically as $kN \to \infty$ for $(k,N) \in 2\Z_{> 0} \times \Sq_r$), and in
fact we know good error estimates.

Consider a Hecke operator $T_p$ ($p \nmid N$) acting on $S_k^\new(N)$ or some 
$S_k^{\new,\eps}(N)$.
Then each Galois orbit $\calO$ of newforms in these spaces corresponds to a factor 
$g_{\calO,p}(x)$ of the characteristic polynomial $c_{T_p}(x) \in \Z[x]$ of $T_p$.
More precisely, we must have $g_{\calO,p}(x) = m_{\calO,p}(x)^j$ for some irreducible polynomial $m_{\calO,p}(x) \in \Z[x]$, and the roots of $g_{\calO,p}(x)$ are in one-to-one 
correspondence with the $T_p$-eigenvalues of the newforms in $\calO$.

Thus a random model for the number and degrees of factors of $c_{T_p}(x)$ will provide a 
 simple heuristic for an upper bound on the number and sizes of Galois orbits of 
 $S_k^\new(N)$.  
 In fact, the main result of \cite{KSW} implies that for a given Galois orbit $\calO$,
 $g_{\calO,p}(x) = m_{\calO,p}(x)$ for 100\% of primes $p$.  Hence we may use
a model for the factorization of $c_{T_p}(x)$ for an arbitrary $p \nmid N$ to model the
Galois orbits of $S_k^\new(N)$.  
 
We recall two heuristic principles on random polynomials, which are quite robust to
the model being considered.

\begin{enumerate}
\item [(RP1)] In the absence of simple reasons for nontrivial factors,
the probability that a well-behaved random polynomial in $\Z[x]$ is irreducible
tends to 1 (typically quickly) as the size of the polynomial grows.

\item [(RP2)] Asymptotically, the probability that a well-behaved random polynomial in 
$\Z[x]$ is reducible over $\Q$ is proportional to the probability that it has a linear
factor over $\Q$.
\end{enumerate}

We do not attempt to define the notions of ``well-behaved'' or the ``size'' of the 
polynomial---indeed we use them in a somewhat vague sense here---but just remark 
that by size we have in mind some combination of size of the 
coefficients and the degree of the polynomial.  The principle (RP1) has been long
studied, and there are many results in this direction.  
See, e.g., \cite{BBBSWW} for a recent study of (RP2),
which the authors term universality.

Apart from the decomposition of $S_k^\new(N)$ into Atkin--Lehner eigenspaces,
there are no obvious reasons why the characteristic polynomials of
 $T_p$'s acting on $S_k^\new(N)$ should factor for $p \nmid N$. 
  Thus (RP1) suggests that,
for $p \nmid N$, $T_p$ acts irreducibly on each Atkin--Lehner eigenspace 100\% of the
time.  This is our first heuristic why \cref{gal-conj} should be true.

To be more precise, we recall a simple model for Hecke polynomials on 
Atkin--Lehner eigenspaces recently proposed by Roberts \cite{roberts}.

Consider the collection $\mathcal P_d(t)$ of degree $d$ monic polynomials
in $\Z[x]$ all roots real size at most $t$. 
It follows from \cite[Theorem 4.1]{AP1} that
$|\mathcal P_d(t)|$ is approximately
\begin{equation} \mathcal R_d(t) := (2t)^{\frac{d(d+1)}2} \prod_{j=1}^d \frac{(j-1)!^2}{(2j-1)!} =2^d t^{\frac{d(d+1)}2} \prod_{j=1}^{d-1} \left(\frac{j}{2j+1} \right)^{d-j}
\end{equation}
when $d+t$ is large (cf.\ \cite[Theorem 3.1]{AP2}).  Thus the probability
that a random polynomial in $\mathcal P_d(t)$ has a factor of degree $e \le \frac d2$ is
approximately 
\begin{multline}
P_{d,e}(t) := \frac{\mathcal R_e(t)  \mathcal R_{d-e}(t)}{2^\delta \mathcal R_d(t)} \\
= \frac 1{2^\delta t^{e(d-e)}} \prod_{j=1}^{e-1} \left( \frac {2j+1}j \right)^{j} 
\prod_{j=e}^{d-e-1} \left( \frac {2j+1}j \right)^{e} 
\prod_{j=d-e}^{d-1} \left( \frac {2j+1}j \right)^{d-j},
\end{multline}
where $\delta = 1$ if $d=2e$ and $0$ otherwise.

Note (for $d > 2$)
\[ P_{d,1}(t) = \frac 1{t^{d-1}} \prod_{j=1}^{d-1} \frac{2j+1}j \]
and
\[ P_{d,e}(t) < \left( \frac 3t \right)^{e(d-e)}. \]
Consequently, we see that both $P_{d,e}(t) \to 0$ and
$P_{d,1}(t) \gg \sum_{e=2}^{\lfloor d/2 \rfloor} P_{d,e}(t)$ as $d+t \to \infty$,
provided $t > 3$.  In fact, by refining the above bound, this analysis works for
fixed $t > 2$ and $d \to \infty$.  Thus (RP1) and (RP2) hold for this model with 
fixed $t > 2$ and $d \to \infty$.

Let $n_\eps = \dim S_k^{\new,\eps}(N) \approx \frac {(k-1)\phi(N)}{12 \cdot 2^r}$.
Then a crude model for the characteristic polynomial
of $T_p$ acting on $S_k^{\new,\eps}(N)$ is a random polynomial in
$\mathcal P_{n_\eps}(2p^{(k-1)/2})$.  One obvious defect is that, even for
fixed $N$ and $\eps$, the probability of reducibility decreases rapidly as 
$p \to \infty$.  

For odd $N$, we will just use a model for the factorization of
$c_{T_2}$ as a model for the Galois orbits.  In light of both \cite{KSW} and
observed data (see \cref{tab6}), the factorization of $c_{T_2}$ does seems to be a very good model
for the sizes of Galois orbits.  Further, as in \cite{roberts}, we may view 
$P_{n_\eps,e}(2^{\frac{k+1}2})$ as a rough model for the probability that 
$S_k^{\new,\eps}(N)$ has a Galois orbit of  size $e$ (or a collection of smaller Galois orbits whose sizes sum to $e$), even when $N$ is even.

This model is not very accurate---as pointed out in
\cite{roberts}, it severely underpredicts the actual number of factorizations of
$c_{T_p}$'s.  In fact it suggests all but finitely many Atkin--Lehner eigenspaces
consist of a single Galois orbit as $N \to \infty$ along $\Sq_r$, which should
not be true at least in weight 2  (cf.\ \cref{sec:ec} and \cref{q:fin}).
We will partially address this by also considering arithmetic statistics in
\cref{sec:ec} below.
However, the model at least suggests the following principles which we believe
in accordance with the data and general expectations about randomness.  
In the following statements, we consider
$r \ge 1$ fixed, $N \in \Sq_r$, $k \ge 2$ even, $\eps$ a sign pattern
for $N$ and $p \nmid N$.  

\begin{enumerate}

\item [(wt)] Given $N$, $\eps$ and $p$, the probability that $T_p$ acts reducibly on $S_k^{\new,\eps}(N)$ decreases rapidly as $k$ becomes large.

\item [(lev)] Given $\eps$ and $p$, the probability
that $T_p$ acts reducibly on $S_k^{\new,\eps}(N)$ decreases rapidly as $N$ grows
coprime to $p$.

\item [(lin)] If $k+N$ is large, the probability that $T_p$ acts reducibly on
$S_k^{\new,\eps}(N)$ is roughly equal to the probability that 
$T_p$ acting on $S_k^{\new,\eps}(N)$ has a rational eigenvalue. 
\end{enumerate}

Note (wt) and (lev) follow from (RP1) since increasing $k+N$ increases
the dimensions of the Atkin--Lehner subspaces.  
Similarly, (lin) follows from (RP2).
Specifically, (lev) suggests \cref{gal-conj}.

Before we discuss other heuristics and data, we discuss the relation of the above
heuristics with the following recent conjecture of \cite{roberts} restricted to our setting
of squarefree level. (In particular, we do not need to account
for quadratic twist classes or CM forms.)  By a rational newform, we mean a newform in $S_k(N)$
with rational Fourier coefficients, i.e., a newform whose Galois orbit has size 1.  
Let $\mathcal F_k$ be the set of
newforms which lie in some $S_k(N)$ with $N$ squarefree.

\begin{conj}[Roberts] Fix $k \ge 6$.  There are only finitely many rational newforms in 
$\mathcal F_k$.  Further, if  $k \ge 52$, there are no rational newforms in $\mathcal F_k$.
\end{conj}

Roberts' support for his conjecture comes from his rough heuristic model together
with an apparent lack of motivic sources for rational newforms in higher weight and
computations in weights $k \le 50$ and levels $N \le C_k$.
Here the bound $C_k$ depends on $k$---e.g., $C_6 = 1000$, $C_{10} = 450$,
$C_{20} = 150$, $C_{30} = 100$ and $C_{40} = 30$.

Of course when $k=2$, rational newforms correspond to isogeny classes of elliptic
curves, so we expect infinitely many rational newforms (cf.\ \cref{sec:ec}).
For $k=4$, Roberts remarks that it is unclear if there should be infinitely many rational 
newforms or not, and that this is related to the existence of suitable Calabi--Yau
threefolds.

Let $\omega(N)$ denote the number of prime divisors of $N$.
Then (lin) suggests the following more speculative conjecture, which implies something 
much stronger than \cref{gal-conj} if Roberts' conjecture is true for $k$.

\begin{conj} \label{conj:finite}
Fix $k$.  Suppose there are only finitely many rational newforms in $\mathcal F_k$.
Then $\# \Orb(S_k^{\new}(N)) = 2^{\omega(N)}$ for all but finitely many squarefree $N$.
\end{conj}

As some numerical evidence for this conjecture,
Roberts observed that in his calculations for $k \ge 6$ that there were only 
four squarefree levels where an Atkin--Lehner eigenspace has multiple Galois orbits 
with no orbits of size 1, which is in line with (lin).

Moreover, in light of the second part of Roberts' conjecture, the above heuristics suggest
that for $k$ sufficiently large, each Atkin--Lehner eigenspace may only be a single Galois orbit for arbitrary squarefree $N$.  Put another way, it is possible
that the generalized Maeda conjecture for squarefree level 
is true with a uniform bound on the
weight: there exists some absolute $k_0$ such that for any $r \ge 0$,
$\#\Orb(S_k^\new(N)) = 2^{\omega(N)}$ for all $k \ge k_0$ and all squarefree $N$.  (The conjecture in \cite{tsaknias} merely asserts that this should be true
if we allow $k_0$ to depend upon $N$.)

\subsection{Elliptic curves} \label{sec:ec}
As mentioned before, random polynomial models as above seem too crude 
to accurately  predict the frequency of small Galois orbits.
One perspective is that the existence of small Galois orbits is due to
the existence of suitable motives, which seem to be hard to model without a
deep understanding of arithmetic.  However, we can supplement 
the random polynomial model above with heuristics from arithmetic geometry.

In particular, based on the principle (wt), we expect that \cref{gal-conj} should be true if it
is true in weight 2.  Our data below corroborate this idea, by indicating the average
number of Galois orbits converges to $2^r$ faster the higher the weight is
(see \cref{sec:lmfdb}).
Moreover, by the principle (lin), we expect that \cref{gal-conj} should be
true if, as $N \to \infty$ in $\Sq_r$, 0\% of weight 2 newforms are rational.

By dimension formulas, we know that the number of weight 2 newforms of level
$N \le X$, $N \in \Sq_r$, grows at least as fast as $\frac{X^2}{\log X}$ 
(up to a constant, this is the asymptotic for $r=1$).  However the number of 
isogeny classes of elliptic curves of arbitrary conductor $\le X$ is $O(X^{1+\eps})$
\cite{duke-kowalski}.
(Heuristics of Watkins \cite{watkins} suggest it is actually $O(X^{5/6})$.)
Consequently, 0\% of weight 2 newforms along levels in $\Sq_r$ are rational.
This gives arithmetic support for our belief in \cref{gal-conj}.

In fact, generalizing work of Serre, Binder \cite{binder} showed that, 
for any weight $k$ and fixed degree $A$, 0\% of weight $k$ newforms of levels $N$ 
have a rationality field of degree
$\le A$ for any sequence of $N \to \infty$ with $\omega(N)$ bounded.

\subsection{Galois groups} \label{sec:galois}
We have discussed heuristics for \cref{gal-conj}, which is an on-average
analogue of two aspects of Maeda's conjecture in the level aspect: 
the number of Galois orbits and the irreducibility of the action of $T_p$.  
The remaining part of Maeda's conjecture is that the rationality fields of
newforms are of type $S_n$.  Since random polynomials tend to have
Galois groups $S_n$, it is reasonable to expect that \cref{gal-conj}
also holds with the added statement that $c_{T_p}$ has Galois group 
$S_n$.

There are examples of newforms with rationality fields
whose Galois group is not a full symmetric group, even when the whole 
Atkin--Lehner space has only a single Galois orbit (see \cref{sec:32}).  
Thus at most we
can hope for statements of the form: for almost all newforms in a given
sequence, the rationality field has Galois group $S_n$.  Analogous to
questions about finiteness of number of rational newforms, one can ask
if ``almost all'' here can be interpreted to mean ``all but finitely many'', i.e.,
given a sequence of Galois orbits of newforms (in distinct quadratic
twist classes), 
are the Galois groups of type $S_n$ for all but finitely many orbits?
Unfortunately, we do not have precise enough heuristics or abundant
enough data to speculate about this.


\section{Data}
\label{sec:data}


\subsection{LMFDB data} \label{sec:lmfdb}
First we present some numerical evidence for \cref{gal-conj}
using data from LMFDB \cite{lmfdb}.  LMFDB contains data
for the newforms in $S_k^\new(N)$ whenever $Nk^2 \le 40000$.
Using this data, we computed the average number $A_{k,r}(X)$ of Galois orbits
over the spaces $S_k^\new(N)$ where $N \in \Sq_r$ with $N < X$
for several values of $X \le 10000$, $2 \le k \le 12$ and $1 \le r \le 3$.
The data are summarized in \cref{tab1,tab2,tab3}, corresponding respectively to $r=1$,
2 and 3.  Blank spaces in the tables denote the situations where there is not
sufficient data in LMFDB to compute the averages.

\begin{table}
\caption{Average numbers $A_{k,1}(X)$ of Galois orbits for $r=1$}
\begin{tabular}{r|rrrrrr}
& $X = 250$ &$500$ & $1000$ & $2500$ & $5000$ & $10000$\\
\hline
$k=2$& 2.038 & 2.484 & 2.679 & 2.684 & 2.577 & 2.483\\
4 & 2.057 & 2.042 & 2.030 & 2.016 &  & \\
6 & 1.981 & 2.000 & 2.000 &  &  & \\
8 & 2.000 & 2.000 &  &  &  & \\
10 & 1.981 &  &  &  &  & \\
12 & 1.943 &  &  &  &  & \\
\end{tabular} \label{tab1}
\end{table}

\begin{table}
\caption{Average numbers $A_{k,2}(X)$ of Galois orbits for $r=2$}
\begin{tabular}{r|rrrrrr}
& $X = 250$ &$500$ & $1000$ & $2500$ & $5000$ & $10000$\\
\hline
$k=2$ & 3.243 & 4.386 & 5.292 & 5.615 & 5.608 & 5.442\\
4 & 4.405 & 4.352 & 4.250 & 4.135 &  & \\
6 & 4.108 & 4.069 & 4.042 &  &  & \\
8 & 3.973 & 3.986 &  &  &  & \\
10 & 4.013 &  &  &  &  & \\
12 & 4.000 &  &  &  &  & \\
\end{tabular} \label{tab2}
\end{table}

\begin{table}
\caption{Average numbers $A_{k,3}(X)$ of Galois orbits for $r=3$}
\begin{tabular}{r|rrrrrr}
& $X = 250$ &$500$ & $1000$ & $2500$ & $5000$ & $10000$\\
\hline
$k=2$ & 3.708 & 5.885 & 8.652 & 11.34 & 12.30 & 12.29\\
4 & 7.500 & 8.902 & 9.237 & 8.701 &  & \\
6 & 8.167 & 8.557 & 8.348 &  &  & \\
8 & 8.292 & 8.197 &  &  &  & \\
10 & 8.083 &  &  &  &  & \\
12 & 7.958 &  &  &  &  & \\
\end{tabular} \label{tab3}
\end{table}

Our first remark about the data is that, for fixed $k, r$ 
the rough shape of the graph $A_{k,r}(X)$ as a function of $X$
initially increases with $N$ (corresponding
to the range where some Atkin--Lehner spaces are 0) 
and then is essentially decreasing.  In accordance with our
heuristics, it appears that $A_{k,r}(X)$ tends to $2^r$ faster
the larger $k$ is and the smaller $r$ is, as the dimensions of the
Atkin--Lehner eigenspaces are larger in these situations.  

In particular, the data for the case $k=2$ and $r=3$ are not sufficient
to make it numerically apparent whether the average tends to $2^3$.
However, the data on the whole seems to be in support of 
\cref{gal-conj}, and also the notion that \cref{gal-conj} should be true
if it is true for $k=2$.  In addition, it seems reasonable to expect that, 
for given $k$, the distribution of sizes of Galois orbits along a sequence
of Atkin--Lehner eigenspaces depends primarily on dimension of the 
Atkin--Lehner eigenspaces and not to any significant amount on the number
$r$ of prime factors of $N$.  Thus, at least to our mind, we can be confident
about \cref{gal-conj} if we are in the special case of $k=2$ and $r=1$,
which is what we focus on below.


\subsection{Data for $S_2(N)$, $N$ prime} \label{sec:32}
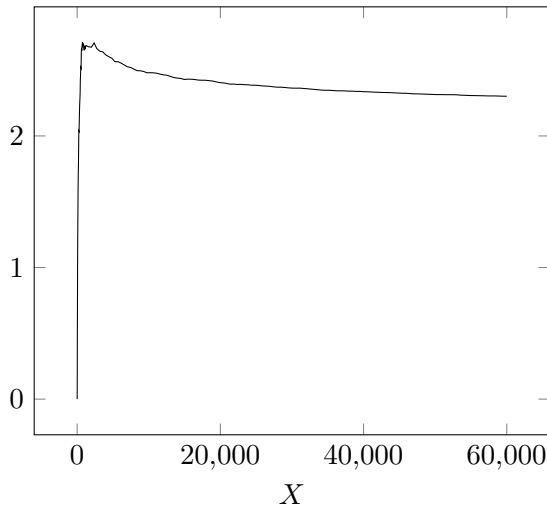
\begin{figure}
\caption{The average number $A_{2,1}(X)$ of Galois orbits for $S_2(N)$, $N < X$ prime}
\begin{tikzpicture}
    \begin{axis}[xlabel=$X$, xtick = {0,20000,40000,60000}, scaled x ticks = false, scaled y ticks = false]
    \addplot[]
        plot coordinates {      
(2, 0.0)
(13, 0.16666666666666666)
(31, 0.5454545454545454)
(53, 0.875)
(73, 1.1904761904761905)
(101, 1.3846153846153846)
(127, 1.5806451612903225)
(151, 1.6944444444444444)
(179, 1.7804878048780488)
(199, 1.8695652173913044)
(233, 2.0392156862745097)
(263, 2.0357142857142856)
(283, 2.081967213114754)
(317, 2.1666666666666665)
(353, 2.23943661971831)
(383, 2.276315789473684)
(419, 2.3333333333333335)
(443, 2.4302325581395348)
(467, 2.4505494505494507)
(503, 2.5208333333333335)
(547, 2.514851485148515)
(577, 2.6320754716981134)
(607, 2.6666666666666665)
(641, 2.646551724137931)
(661, 2.669421487603306)
(701, 2.6904761904761907)
(739, 2.7099236641221376)
(769, 2.6911764705882355)
(811, 2.702127659574468)
(839, 2.6986301369863015)
(877, 2.6887417218543046)
(911, 2.6794871794871793)
(947, 2.670807453416149)
(983, 2.6506024096385543)
(1019, 2.6783625730994154)
(1049, 2.659090909090909)
(1087, 2.679558011049724)
(1109, 2.6774193548387095)
(1153, 2.670157068062827)
(1193, 2.683673469387755)
(1229, 2.6865671641791047)
(1597, 2.6772908366533863)
(1993, 2.6744186046511627)
(2371, 2.7065527065527064)
(2749, 2.663341645885287)
(3187, 2.6430155210643016)
(3581, 2.6387225548902196)
(4001, 2.6152450090744104)
(4421, 2.6023294509151413)
(4861, 2.589861751152074)
(5281, 2.563480741797432)
(5701, 2.5645805592543276)
(6143, 2.553058676654182)
(6577, 2.539365452408931)
(7001, 2.5260821309655936)
(7507, 2.519453207150368)
(7927, 2.5074925074925076)
(8389, 2.4957183634633684)
(8837, 2.4950045413260673)
(9293, 2.4900086880973067)
(9739, 2.479600333055787)
(10181, 2.480415667466027)
(10663, 2.479631053036126)
(11159, 2.4759437453737974)
(11677, 2.4682369735902925)
(12113, 2.46381805651275)
(12569, 2.4610259826782146)
(13009, 2.4526112185686655)
(13513, 2.443472829481574)
(13997, 2.439733494851605)
(14533, 2.4362139917695473)
(14951, 2.4288977727013137)
(15413, 2.4314269850083288)
(15887, 2.4311183144246353)
(16411, 2.42872172540768)
(16921, 2.424397744746284)
(17393, 2.423288355822089)
(17903, 2.42320819112628)
(18329, 2.420752022846264)
(18911, 2.4170153417015343)
(19427, 2.4107223989095865)
(19913, 2.4051532652154597)
(20359, 2.4033029117774882)
(20899, 2.3985538068906846)
(21391, 2.3931695127030403)
(21851, 2.39453284373725)
(22343, 2.391043582566973)
(22853, 2.3920031360250884)
(23327, 2.3894655901576316)
(23831, 2.3889098453413804)
(24317, 2.3857830433172897)
(24889, 2.3853144311159578)
(25409, 2.38307747233131)
(25919, 2.379866713433883)
(26407, 2.378490175801448)
(26893, 2.3771602846492716)
(27457, 2.372875708097301)
(27953, 2.370698131760079)
(28513, 2.37020316027088)
(28949, 2.368454458901936)
(29453, 2.3655107778819118)
(30071, 2.363272839126423)
(30577, 2.3629203271735837)
(31121, 2.36317517159057)
(31607, 2.360776242281682)
(32173, 2.3584468270066647)
(32611, 2.3558983147672095)
(33119, 2.3537031822021968)
(33617, 2.3507359066925853)
(34159, 2.3483976992604765)
(34651, 2.3469332612807348)
(35171, 2.3468408424420155)
(35771, 2.3451723230728754)
(36293, 2.3427681121786548)
(36787, 2.3432453217123816)
(37313, 2.3424449506454064)
(37831, 2.340664833791552)
(38377, 2.3401629227351273)
(38923, 2.338454035601073)
(39439, 2.338713562996868)
(39979, 2.3363484884551298)
(40543, 2.335685721006822)
(41113, 2.3338758428272497)
(41609, 2.3325672259250747)
(42083, 2.331515564644399)
(42569, 2.3304875308919346)
(43063, 2.3299266829593424)
(43669, 2.3291584267194025)
(44203, 2.327320147793958)
(44753, 2.3263814233498175)
(45317, 2.32439906402893)
(45863, 2.32435276783835)
(46451, 2.3220162466152883)
(46997, 2.3199340342197483)
(47533, 2.319730667210773)
(48073, 2.3185215108058976)
(48619, 2.3175364927014597)
(49139, 2.316175014848545)
(49667, 2.3154283473828663)
(50153, 2.314113764317608)
(50767, 2.314170351855412)
(51341, 2.3136545419920016)
(51817, 2.3133371062063763)
(52363, 2.313399364604747)
(52937, 2.312905017589335)
(53453, 2.310218308567235)
(54001, 2.3090347209598256)
(54547, 2.307872455413439)
(55109, 2.3067309409034102)
(55681, 2.3063174659352326)
(56197, 2.3057358358182776)
(56713, 2.305338202051817)
(57193, 2.3040855024995692)
(57751, 2.303879678687404)
(58243, 2.3040162684290797)
(58897, 2.3031423290203326)
(59369, 2.3019496750541575)
(59921, 2.3017683027598745)
(59999, 2.301634472511144)        
        };
\end{axis}
\end{tikzpicture}
 \label{fig1}
\end{figure}

To gather more numerical evidence for \cref{gal-conj}, we computed
$A_{2,1}(X)$ for $X \le 60000$.
These calculations were carried out in parallel using Sage \cite{sage} 
on the University of Oklahoma's supercomputing facilities (OSCER)
over the course of several weeks.
See \cref{fig1} for a graph of $A_{2,1}(X)$, which appears to be
eventually decreasing to 2, as conjectured.
We remark that $A_{2,1}(60000) \approx 2.3016$ (compare with \cref{tab1}).
The apparent slow rate of convergence is expected in light of the well-known
numerical phenomenon that the proportion of weight 2 newforms accounted
for by elliptic curves tends to 0 quite slowly.

One consequence of \cref{gal-conj} would be that, in 100\% of prime levels, $S_2(N)$
has exactly 2 Galois orbits.  (For $N > 59$, $S_2(N)$ has at least 2 orbits.)
In fact, provided the number of Galois orbits
does not grow too fast along any subsequence, this is equivalent to the $k=2$, $r=1$ case
of \cref{gal-conj}.  \cref{tab4} summarizes how often we get exactly 2 (or 3, or 4, etc)
Galois orbits in weight 2 in certain ranges.  
These numerics suggest that indeed there are  exactly 2 Galois orbits 100\% of the time.
We remark that for prime $N < 60000$, the maximum number of Galois orbits is 10.

\begin{table}
\caption{Counts of number of Galois orbits for $S_2(N)$ with $N$ prime}
\begin{tabular}{r| rrrrrr | r}
\multicolumn{1}{c}{} & \multicolumn{6}{c}{{number of orbits}} & \% with  \\
 & $2$ &$3$ & $4$ & $5$ & $6$ & $7+$ & 2 orbits \\
\hline
$0 < N < 10000$ & 777 & 331 & 67 & 25 & 9 & 6 & 63.2\% \\
$10000 < N < 20000$ & 786 & 193 & 39 & 9 & 5 & 1 & 76.1\% \\
$20000 < N < 30000$ & 769 & 176 & 31 & 5 & 1 & 1 & 78.2\% \\
$30000 < N < 40000$ & 768 & 158 & 23 & 7 & 2 & 0 & 80.2\%  \\
$40000 < N < 50000$ & 750 & 162 & 14 & 3 & 1 & 0 & 80.6\%  \\
$50000 < N < 60000$ & 765 & 121 & 28 & 6 & 4 & 0 & 82.8\% \\
\end{tabular}  \label{tab4}
\end{table}

One of our heuristics for \cref{gal-conj} uses the idea (lin), that most of the time an 
Atkin--Lehner eigenspace has multiple Galois orbits, the multiple orbits are
accounted for by the existence of a rational newform.  In \cref{tab5}, we summarize the 
number of small Galois orbits in various ranges, and observe that indeed most of the time 
there is a small Galois orbit, it is of size 1.  In fact, there are no orbits of ``moderate'' size:
for $10000 < N < 60000$, there are no Galois orbits of size $6 \le d \le 300$. 
This suggests the following question:

\begin{question} \label{q:fin}
Fix $k, r, d$.  Are there infinitely many Galois orbits of size $d$
in the union of spaces $S_k^\new(N)$ with $N \in \Sq_r$?
\end{question} 

Note that for $d=1$, this is just asking about the infinitude of rational newforms,
which is the topic of Roberts' conjecture discussed above.  In fact, Roberts' conjecture
together with \cref{conj:finite} would imply that for $k$ large, the answer 
is negative for all $r, d$.
On the other hand, when $k=2$ and $d$ is small, we expect this question has a positive 
answer.  So the most novel case of this question is when $k$ is small but $d$ is not. 
At least for $k=2$, $r=1$ and $d$ sufficiently large, our data suggest the answer 
may be no.

\begin{table}
\caption{Counts of small Galois orbits for $S_2(N)$ with $N$ prime}
\begin{tabular}{r| rrrrrr r}
& \multicolumn{7}{c}{size of orbits} \\
  & $1$ &$2$ & $3$ & $4$ & $5$ & $6$ & $7$ \\
\hline
$0 < N < 10000$ & 329 & 212 & 76 & 28 & 20 & 11 & 18 \\
$10000 < N < 20000$ & 200 & 104 & 16 & 3 & 0 & 0 & 0 \\
$20000 < N < 30000$ & 176 & 80 & 5 & 2 & 1 & 0 & 0 \\
$30000 < N < 40000$ & 171 & 56 & 5 & 1 & 0 & 0 & 0 \\
$40000 < N < 50000$ & 140 & 56 & 7 & 0 & 0 & 0 & 0 \\
$50000 < N < 60000$ & 152 & 57 & 2 & 0 & 0 & 0 & 0
\end{tabular} \label{tab5}
\end{table}

\subsection{The method}
Now we describe our method to compute Galois orbits.  For an odd prime $N$,
let $B = B_N$ be the definite quaternion algebra of discriminant $N$.  Then we computed
the Brandt matrix $T_2^B$ for a maximal order $\calO_B$ of $B$, which acts on the 
space of $M$ quaternionic modular forms associated to $\calO_B$.  
This space of quaternionic
modular forms is Hecke isomorphic to $M_2(N) \simeq \C E_{2,N} \oplus S_2(N)$,
where $E_{2,N}$ is the normalized weight 2 level $N$ holomorphic Eisenstein series.
The Eisenstein eigenvalue of $T^B_2$ is 3, and thus the eigenvalues of $T^B_2$
acting on $B_N$ are 3 together with the eigenvalues of $T_2$ acting on $S_2(N)$.
We compute the characteristic polynomial $c_{T^B_2}(x) = (x-3) c_{T_2}(x)$.

If $c_{T_2}(x)$ has no repeated factors, then the number of Galois orbits in $S_2(N)$
is simply the number of irreducible factors of $c_{T_2}(x)$.
If $c_{T_2}(x)$ has repeated factors, then we repeat the above calculation with
successive $T_p$'s until this method succeeds.  We performed these calculations
in parallel by treating different $N$ on different cores.
Most of the calculation time is spent computing the Brandt matrices $T_p^B$
and their characteristic polynomials, and the computational complexity 
increases both with $N$ and with $p$.  
For $N$ close to 60000, this calculation for a single $T_p$ 
can take over 24 hours of wall time.

\subsection{Irreducibility of Hecke polynomials} \label{sec:irred}
In most cases, $T_2$ acts on $S_2(N)$ with no repeated eigenvalues.  
Even when $T_2$ does not, we typically do not
have to try many $T_p$'s to find one that does.  \cref{tab6} shows for how many
primes $N < 60000$ a given $p$ is minimal such that $T_p$ has no
repeated eigenvalues.

\begin{table}
\caption{Frequency that $p$ is the smallest prime such that $T_p$ acts on $S_2(N)$ with no repeated eigenvalues ($N < 60000$ prime)}
\begin{tabular}{c|rrrrrrrrrr}
$p$ &  2 & 3 & 5 & 7 & 11 & 13 & 17 & 19 & 23 & 47 \\
 \hline
frequency & 5815 & 158 & 42 & 14 & 15 & 2 & 3 & 4 & 3 & 1  
\end{tabular} \label{tab6}
\end{table}

More generally, given $p$ we can ask how many prime $N \ne p$
are there such that $c_{T_p}$ has no repeated roots for $S_2(N)$?
Note that if we have multiple rational newforms occurring in $S_2(N)$, the naive 
probability that two given such newforms $f_1$ and $f_2$ have the same $a_p$ is 
approximately $\frac 1{4 \sqrt p}$ by Deligne's bound.  Since we expect that there are
infinitely many prime levels $N$ where $S_2(N)$ has 2 rational newforms (e.g.,
coming from Neumann--Setzer elliptic curves \cite{setzer}) we expect that $c_{T_p}$
has repeated roots for infinitely many $N$.  

To avoid this situation, let us restrict to the case where $S_2(N)$ has exactly 2 Galois
orbits.  Out of the 4615 prime levels $N < 60000$
such that $S_2(N)$ has exactly 2 Galois orbits, there is only one level $N$
such that $T_2$ has repeated eigenvalues on $S_2(N)$, namely $N=251$.
Here $T_2$ acts irreducibly on the 17-dimensional root number $+1$
subspace of $S_2(251)$ but acts reducibly on the 4-dimensional root number $-1$
subspace.  (Incidentally, a newform in the latter space has rationality field with Galois group $D_8$.)  Based on the rarity of repeated eigenvalues of $T_2$, 
we guess that there may be no other such $N$.

This is related to a question studied in \cite{murty} and \cite{KSW}: 
given a non-CM newform $f \in S_k(N)$ 
with rationality field $K$, how often is $\Q(a_p(f))$ a proper subfield of $K$?
Restricted to our setting of squarefree level and trivial nebentypus,
 \cite[Conjecture 3.4]{murty} asserts that this happens for infinitely many $p$ exactly
in the following cases: $k=2$ and $K$ is quadratic, cubic or quartic with a quadratic
subfield.  See also \cite{vanhirtum} for the case $k=[K:\Q] = 2$.
In the cases where we expect $\Q(a_p(f)) = K$ for all but finitely many $p$,
we can ask if something stronger is true.

\begin{question} \label{q:ratl}
Given $k$ and squarefree $N$, let $f \in S_k(N)$ be a non-CM newform with rationality field $K$.
Suppose $k \ge 4$ or $[K:\Q] \ge 5$.
Is $\Q(a_p(f)) = K$ for all $p \nmid N$?
\end{question}

Again, the answer may depend on $k$ and $K$, and possibly $N$,
but it is not entirely unreasonable to think this question may have a positive answer 
if some combination $k$, $K$ and $\Gal(K/\Q)$ are sufficiently large.
At least when $N=1$, there is some evidence towards a positive answer
(e.g., \cite{VX}, \cite{bengoechea}).  Note that an affirmative answer to 
\cref{q:ratl} under appropriate
conditions would imply that $T_p$ acts irreducibly on each Galois orbit
of $S_k(N)$ for all but finitely many squarefree $N$ coprime to $p$.


\section{Proof of \cref{thm1}}
\label{sec:Lvals}

Fix $k \ge 2$ even.
For a prime level $N$, the Atkin--Lehner eigenspaces of $S_k^\new(N)$ are
simply the spaces with fixed root number $\pm 1$.

Let $\sigma \in \Gal(\bar \Q/\Q)$.
It follows from an algebraicity result of Shimura \cite[Theorem 1]{shimura:1976}
that $L(k, f^\sigma) \ne 0$ if and only if $L(k,f) \ne 0$.
Thanks to the extension of the Gross--Zagier formula by 
Zhang, we also know that $L'(k, f^\sigma) \ne 0$ if and only if
$L'(k,f) \ne 0$ \cite[Corollary 0.3.5]{zhang}.  

Consequently to show \cref{thm1}, it suffices to know that for $N$ sufficiently large,
there exist $f \in S_k(N)$ with $L(k,f) \ne 0$ (so $f$ necessarily has root number
$+1$) and $g \in S_k(N)$ with $L'(k,g) \ne 0$ and $g$ having root number $-1$.
This follows, e.g.,  from the works \cite{vanderkam}, \cite{kowalski-michel}
and \cite{iwaniec-sarnak} mentioned in the introduction.


%
%

\end{document}